\documentclass[12pt]{article}
\newtheorem{thm}{Theorem}
\newtheorem{lem}[thm]{Lemma}
\newtheorem{prop}[thm]{Proposition}
\newtheorem{cor}[thm]{Corollary}
\newcommand{\Aut}{\mathrm{Aut}}
\newcommand{\Int}{\mathrm{Int}}
\newcommand{\Ad}{\mathrm{Ad}}
\newcommand{\mod}{\mathrm{mod}}
\newcommand{\Ker}{\mathrm{Ker}}
\newcommand{\id}{\mathrm{id}}
\renewcommand{\epsilon}{\varepsilon}
\usepackage{latexsym}
\author{MASUDA Toshihiko \\
Graduate School of Mathematics, Kyushu University, \\
Fukuoka, 812-5250, JAPAN}
\date{}
\title{Evans-Kishimoto type argument for actions of discrete amenable
groups on McDuff factors} 
\begin{document}
\maketitle
\begin{abstract}
We apply the Evans-Kishimoto type argument for centrally free actions of
 discrete amenable  groups on McDuff factors, and classify them.
 Especially, we present a different proof that 
 the Connes-Takesaki
 modules are complete cocycle conjugacy invariants for centrally free actions of
 discrete amenable groups on injective factors. 
\end{abstract}
\section{Introduction}
In the theory of operator algebras, the study of automorphism groups is
one of the most important subjects. Especially, since Connes succeeded in classifying
automorphisms of the approximately finite dimensional (AFD) factor of
type II$_1$ in \cite{Co-peri} and \cite{Con-auto}, 
classification of actions of discrete amenable groups on injective
factors has been solved in 
\cite{J-act}, \cite{Ocn-act}, \cite{Su-Tak-act}, \cite{KwST} and
finally in \cite{KtST}. 

The idea of Connes's classification is following. 
First he constructed
tensor product type model automorphisms (or actions) on the AFD type
II$_1$ factor. Then he showed ``the model action splitting'', i.e., 
every automorphism contains model
actions as tensor product components after an appropriate inner 
perturbation, 
and then proved that it is cocycle conjugate to the model action. 
In his arguments, the main tool is a Rohlin property for automorphisms,
(or actions). In \cite{Con-auto}, he showed the noncommutative
version of Rohlin type theorem for a certain class of automorphisms. By
means of the Rohlin type theorem,
he proved the stability (or 1-cohomology vanishing theorem) for
automorphisms, which is the important step 
in his classification. 
Connes's argument has been developed by Jones for finite group case in \cite{J-act}, 
and Ocneanu for general discrete amenable group case in \cite{Ocn-act}. 
Especially, the Rohlin type theorem was extended to the case of
discrete amenable groups by Ocneanu, and he proved
several cohomology vanishing theorems in \cite{Ocn-act}. 

On the other hand, another approach has been made in the study of
automorphisms of $C^*$-algebras. In \cite{EvKi}, Evans and Kishimoto 
developed the intertwining argument for classification of automorphisms
with the Rohlin property. In their approach, they compare given two
automorphisms directly without using model actions. 
As a consequence,
they can treat a wide class of automorphisms, and obtained classification
results. (In $C^*$-algebra case, arguments based on  model actions 
forces us to make a strict restriction on actions.)
Their intertwining argument has been farther developed in \cite{Naka} for
automorphisms of purely infinite simple $C^*$-algebras, and also for
finite group actions in \cite{Iz-roh1}. 

In this paper, we apply the Evans-Kishimoto type intertwining argument
for actions of discrete amenable groups on McDuff factors based on Ocneanu's
Rohlin type theorem.  
Our main theorem says that if two centrally free actions of a discrete
amenable group on a McDuff factor differ up to approximately inner
automorphisms, then they are cocycle conjugate.
Especially if we apply our result to injective factors, then we get
the complete classification of centrally free actions of discrete
amenable groups in terms of Connes-Takesaki module by using  
the characterization of approximately inner automorphisms
for injective factors in \cite{KwST}. 
Hence this is an another proof of classification results in \cite{Con-auto},
\cite{J-act}, \cite{Ocn-act}, \cite{Su-Tak-act}, \cite{KwST} 
for centrally free actions. However our approach seems to be more  
unified and simple one, and this is an advantage of our theory. 

Our result is also applicable to the
classification of group actions on subfactors by suitable modification. 
For examples, we present a different proof of Popa's result in 
\cite[Theorem 3.1]{Po-act}. (We remark that the classification result of strongly
amenable subfactors of type II$_1$ by Popa in \cite{Po-act} is crucial in 
our argument.)

\textbf{Acknowledgements.} The author is grateful to Professor Izumi, 
Professor Katayama, and Professor Takesaki for comments on this work. He
is supported by Grant-Aid for Scientific Research, Japan Society for the
Promotion of Science. 

\section{Preliminaries and notations}
Let $M$ be a von Neumann algebra.
For $\varphi\in M_*^+$, we let $\|x\|_\varphi=\sqrt{\varphi(x^*x)}$,
$\|x\|_\varphi^{\#}= \sqrt{(\|x\|_\varphi^2+\|x^*\|^2_\varphi)/2}$,
$|x|_\varphi=\varphi(|x|)$. Note $|x|_\varphi$ is not necessary a norm
unless $\varphi$ is tracial, since it is not subadditive.  
For $x\in M$ and $\varphi\in M_*$, functionals $x\varphi$ and $\varphi x$
are defined as $x\varphi(y)=\varphi(yx)$ and $\varphi
x(y)=\varphi(xy)$ respectively. We let $[x, \varphi]=x \varphi - \varphi x $.
To avoid possible confusions, we also denote $x\varphi $ and $\varphi x$ 
by $x\cdot \varphi$ and $\varphi\cdot x$. We denote
$\phi\circ\alpha^{-1}$ by 
$\alpha(\phi)$ for $\phi\in M_*$ and $\alpha\in \Aut(M)$.

We use the notation $A\subset \subset B$ if $A$ is a finite subset
of $B$, and denote the cardinality of $A$ by $|A|$.

Fix a free ultrafilter $\omega$ over $\mathbf{N}$.
We define $M_\omega$ and $M^\omega$ as in \cite{Ocn-act}. Each $\alpha
\in \Aut(M)$ gives an automorphism $\alpha^\omega\in
\Aut(M^\omega)$, and $\alpha^\omega\mid_{M_\omega}\in \Aut(M_\omega)$. 
For $\varphi\in M_*$ and $X=(x_n)\in M^\omega$,
$\varphi^\omega(X):=\lim_{n\rightarrow 
\omega}\varphi(x_n)$ is a normal functional of $M^\omega$, which we
denote by $\varphi$ for simplicity. 
When $M$ is a factor,
$\tau_\omega(X):=\lim_{n\rightarrow
\omega}x_n(\in \mathbf{C})$ always exists in $\sigma$-weak topology
for $X=(x_n)\in M_\omega$, and
$\tau_\omega$ is a tracial state on $M_\omega$.
We denote by $|X|_1$
the $L^1$-norm with respect to this trace.

We collect
fundamental,  frequently used
 inequalities in this paper.
\begin{lem}\label{lem:ineq}
 We have the followings for $\varphi\in M_*^+$, $x\in M$. \\
$(1)$ $\|x\cdot\varphi\|\leq\sqrt{\|\varphi\|} \|x\|_\varphi$, 
$\|\varphi\cdot x\|\leq \sqrt{\|\varphi\|}\|x^*\|_\varphi$, 
$\|[x,\varphi]\|\leq 2\sqrt{\|\varphi\|}\|x\|_\varphi^{\#}$. \\
$(2)$ $\|x\|_\varphi^2\leq \|x\cdot \varphi\|\|x\|$, 
$\|x^*\|_\varphi^2\leq \|\varphi\cdot x\|\|x\|$.\\ 
$(3)$ $\|x\|_\varphi^{\#}\leq \sqrt{\frac{1}{2}(|x|_\varphi+|x^*|_\varphi)\|x\|}$. \\
$(4)$ For $x_i \in M^\omega$ and $y_i\in M_\omega$, 
$|\sum_ix_iy_i|_\varphi\leq \sum_i\|x_i\||y_i|_1$.
\end{lem}
\textbf{Proof.} It is elementary to see (1), (2), (3). 
See \cite[Lemma 7.1]{Ocn-act} for the proof of (4).
\hfill$\Box$

Next we recall Ocneanu's Rohlin type theorem, which is a main tool in
the proof of Lemma \ref{lem:stab} below.

\begin{thm}[{\cite[Theorem 6.1]{Ocn-act}}] 
  Let $M$ be a McDuff factor, 
$G$ a discrete amenable group, $\alpha$ an action of $G$ on $M_\omega$
 which is strongly free and semiliftable. Let $\varepsilon>0$ and let
$\{K_i\}_{i\in I}$ be an
 $\varepsilon$-paving  family. Then there exists a partition of unity
 $\{E_{i,k}\}_{i\in I, k,l\in K_i} \subset M_\omega$ such that 
$$\sum_{i\in I}|K_i|^{-1}\sum_{k\in K_i}|\alpha_{kl^{-1}}(E_{i,l})-E_{i,k}|_1\leq 
5\varepsilon^\frac{1}{2},$$
$$[\alpha_g(E_{i,k}), E_{j,l}]=0, \mbox{ for all } g\in G, i,j\in I, k\in K_i,l\in K_j.$$  
\end{thm}

See \cite{Ocn-act} for terms appearing in the above theorem. Here we
briefly explain how to construct an $\varepsilon$-paving family $\{K_i\}$. Fix $N\in
\mathbf{N}$ such that $N >\frac{4}{\varepsilon}\log \varepsilon^{-1}$,
and  set $\delta:=(\varepsilon/3)^N$.  Let $K_{n+1}$ be a 
$(\delta|\bar{K}_n|^{-1}, \bar{K}_n)$ invariant finite set, where
$\bar{K}_n:=\cup_{1\leq i\leq n} K_{i}$. 
(In this paper, we say $K$ is $(\varepsilon, F)$-invariant if $|K\cap
\bigcap_{g\in F}g^{-1}K|\geq (1-\varepsilon)|K|$.)
Then $\{K_i\}_{1\leq i\leq N}$
is shown to be an $\varepsilon$-paving family. In this construction, 
each $K_i$ can be chosen arbitrarily invariant.

\section{Classification}

We state the main result in this paper.

\begin{thm}\label{thm:main}
Let $M$ be a McDuff factor, $G$  a countable discrete amenable group,
 $\alpha$, $\beta$ centrally free actions of $G$ on $M$. 
 If $\alpha_g\beta^{-1}_g\in \overline{\Int}(M)$ for every $g\in G$,
 then there exist an $\alpha$-cocycle $v_g$ and $\theta\in \overline{\Int}(M)$ such that  
 $\Ad v_g\alpha_g=\theta\circ\beta_g\circ \theta^{-1}$. Moreover 
 we can choose $v_g$ close to $1$, i.e., for given any
 $\varepsilon>0$, $F\subset\subset G$ and $\varphi\in M_*^+$,
we can choose $v_g$ so that
 $\|v_g-1\|_\varphi< \varepsilon$ for any $g\in F$.
\end{thm}

The assumption $\alpha_g\beta_g^{-1}\in \overline{\Int}(M)$ implies that we can
approximate $\alpha_g$ by the inner perturbation of $\beta_g$. However
the inner perturbation of $\beta_g$ is not necessary an action of
$G$. Hence our first task is to approximate $\alpha$ by the cocycle
perturbation of $\beta$.

\begin{prop} Under the assumption in Theorem \ref{thm:main},
 there exist $\beta$-cocycles $u^n_g$, $n=1,2,\cdots$ such that
 $\alpha_g=\lim_{n\rightarrow \omega}\Ad u^n_g \beta_g$.
\end{prop}
\textbf{Proof.} Since $\alpha_g\beta_g^{-1}\in \overline{\Int}(M)$,
there exist unitaries $u_g^n$, $g\in 
G$, $n=1,2,\cdots$ such that $\alpha_g=\lim_{n\rightarrow\infty}\Ad
u_g^n\beta_g$. Set $U_g:=(u^n_g)\in M^\omega$. Then $\alpha_g^\omega=\Ad
U_g\beta^\omega_g$ holds on $M(\subset M^\omega)$. Set
$u(g,h)=U_g\beta_g^\omega(U_h)U_{gh}^*$. Then it is easy to check $u(g,h)\in M_\omega$ and 
$\{\Ad U_g\beta_g^\omega|_{M_\omega}, u(g,h) \}$ is a cocycle action on
$M_\omega$. Moreover $\Ad U_g\beta_g^\omega$ is strongly free in the
sense of \cite[Definition 5.6]{Ocn-act} by
\cite[Lemma 5.7]{Ocn-act}. 
By Ocneanu's 2-cohomology vanishing theorem
\cite[Proposition 7.4]{Ocn-act},
we get $c_g\in 
U(M_\omega)$ such that $c_g\Ad U_g\beta_g^\omega(c_h)u(g,h)c_{gh}^*=1 $. 
Let $c_g=(c_g^n)$ be a representing sequence consisting of unitaries.
Set $W_g:=c_gU_g=(c_g^nu_g^n)$. 
Since $(c_g^n)$ is a centralising
sequence, $\Ad\, c_g^n$ converges to $\id_M$. Hence
$\alpha_g=\lim_{n\rightarrow\omega}\Ad c_g^nu_g^n\beta_g$,  
and $W_g$ is a $\beta^\omega$-cocycle. Set
$w^n_g=c_g^nu_g^n$ and 
$w'_n(g,h):=w_g^n\beta_g(w^n_h)w^{n*}_{gh}$. 
Then $(\Ad w_g^n\beta_g, w_n'(g,h))$ is a 
cocycle action, and the cocycle identity $W_g\beta_g(W_h)=W_{gh}$
implies  $\lim_{n\rightarrow \omega}w'_n(g,h)=1$ in the $\sigma$-strong*
topology.  
By Ocneanu's
2-cohomology vanishing theorem \cite[Theorem 7.6]{Ocn-act}, we can find
$d_g^n\in U(M) $ such that  
$d_{g}^n\Ad w_g^n\beta_g(d^n_h)w'_n(g,h)d_{gh}^{n*}=1$ and
$\lim_{n\rightarrow\omega }d_g^n=1$ in the $\sigma$-strong*
topology. Then $d_g^nw_g^n$ is a  
$\beta_g$-cocycle and $(d^n_gw^n_g)=(w^n_g)$ in $M^\omega$. 
It is easy to see $\alpha_g\beta_g^{-1}=\lim_{n\rightarrow \omega}\Ad
d_g^nw^n_g$. 
\hfill$\Box$.

We get the following corollary immediately. 
\begin{cor}
 Let $\alpha$, $\beta$ be as in Theorem \ref{thm:main}. Then for any 
 $\varepsilon>0$, $\Phi\subset\subset M_*$ and $F\subset\subset G$,
 there exists a  $\beta$-cocycle $u_g$ such 
 that $\|\alpha_g(\phi)-\Ad u_g\beta_g(\phi)\|<\varepsilon $ for every $g\in F$
 and $\phi\in \Phi$. 
\end{cor}

Next we  show an almost 1-cohomology vanishing, which
is essential in our argument. 
\begin{lem}\label{lem:stab} Let $\alpha$ be a centrally free action of
 $G$ on $M$.
 For any $\varepsilon>0$,  $F=F^{-1}\subset\subset G$, 
$\Phi^+\subset\subset M_*^+$, 
and $\Phi \subset\subset M_*$, there exist 
 $\delta>0$ and $\Psi\subset\subset M_*$ with the following property;
 for any $\alpha$-cocycle $\{v_g\}$ with 
$\|[v_g,\psi]\|<\delta$, $g\in F$, $\psi\in \Psi $, we
 can find 
 $w\in U(M)$ such that $\|[w,\varphi]\|<\varepsilon $ for every $\varphi
 \in \Phi$, 
 and $\|v_g\alpha_g(w^*)w-1\|_{\phi}^{\#}<\varepsilon$ for every $g\in F$ 
 and $\phi\in \Phi^+$. When $\Phi=\emptyset$, then $\Psi=\emptyset$ is possible.  
\end{lem}
\textbf{Proof.} 
Since every $\phi\in M_*$ is decomposed as
$\phi=\phi_1-\phi_2+i(\phi_3-\phi_4)$, $\phi_i\in M_*^+$, it is
sufficient to show the lemma in the case
$\Phi=\Phi^+\subset\subset M_*^+$. 
We may assume $0<\epsilon <1$ and $\|\phi\|\leq 1$, $\phi\in \Phi^+$. 
Fix $\epsilon'>0$ with $\epsilon'<(\epsilon/8)^4$.
Let $\{K_i\}_{i\in I}$ be an $\epsilon'$-paving family such that each $K_i$
is $(\epsilon', F)$ invariant. 
We may assume that $K_i$'s are
in a subgroup of $G$ generated by $F$. 
Define $\mathrm{Length}(g):=\min\{n\mid g=h_1h_2\cdots h_n, h_i\in F\}$, 
and define $L$ by
$$L:=\max\{\mathrm{Length}(g)\mid g\in K_i,  i\in I\}.$$ 
Fix $\delta> 0$ so that $\sum_i|K_i|(L+1)\delta<\epsilon/3$. Define
$\Psi$ by 
$$\Psi:=\Phi^+\cup\bigcup_{1\leq k\leq L-1, \atop{g_i\in F}
}\alpha_{g_1g_2\cdots g_k}^{-1}(\Phi^+).$$ 
By Ocneanu's Rohlin Theorem 
there exists a partition of
unity $\{E_{i,k}\}_{i\in I, k\in K_i}\subset M_\omega$ such that 
$$\sum_{i\in I}|K_i|^{-1}\sum_{k,l\in K_i}|\alpha_{kl^{-1}}(E_{i,l})-E_{i,k}|_1
<5\epsilon^{'\frac{1}{2}},$$  
$$ [\alpha_{g}(E_{i,k}), E_{j,l}]=0\,\, \mbox{for all }g\in G,
 i,j\in I, k\in K_i, l\in K_j. 
$$
Then by \cite[Corollary 6.1]{Ocn-act}, we have followings.
$$\sum_{i\in I}\sum_{k\in K_i\cap g^{-1}K_i}|\alpha(E_{i,k})-E_{i,gk}|_1
\leq 10\epsilon^{'\frac{1}{2}}\mbox{ for every } g \in F,
$$
$$\sum_{i\in I}\sum_{k\in A_i}|E_{i,k}|_1\leq
\delta'+5\epsilon^{'\frac{1}{2}},$$ 
for any subsets $A_i\subset K_i$ with $|A_i|\leq \delta'|K_i|$, $i\in I$.

Let $v_g$ be an $\alpha$-cocycle with $\|[v_g,\phi]\|<\delta$,
$g\in F$, $\phi\in \Psi$.
Define $W:=\sum_{i,k}v_k^*E_{i,k}$. 
First we estimate $\|v_g\alpha_g(W^*)W-1\|_\phi^\#$ as in the proof of
\cite[Proposition 7.2]{Ocn-act}. 
We investigate $|v_g\alpha_g(W^*)W-1|_\phi$ and 
$|(v_g\alpha_g(W^*)W-1)^*|_\phi$, then use Lemma \ref{lem:ineq}.

We divide
$v_g\alpha_g(W^*)W-1=\sum_{i\in I,k\in K_i}\sum_{j\in I,l\in K_j}(v_g\alpha_g(v_k)v_l^*-1) 
\alpha_g(E_{i,k})E_{j,l}$ into three parts as follows. 
\begin{eqnarray*}
\sum_{i\in I,k\in K_i}\sum_{j\in I,l\in K_j}(*)&=&\sum_{j\in I,l\in K_j}
\sum_{i,k\in K_i\cap g^{-1}K_i}(*)+
 \sum_{j\in I,l\in K_j}\sum_{i\in I,k\in K_i\backslash g^{-1}K_i}(*) \\
&=& \sum_{j\in I,l\in K_j}\sum_{j\ne i\in I, \atop{k\in K_i\cap g^{-1}K_i}}(*)+
    \sum_{j\in I,l\in K_j,\atop{k\in K_j\cap g^{-1}K_j}}(*)
 +\sum_{j\in I,l\in K_j}\sum_{i,k\in K_i\backslash g^{-1}K_i}(*) \\
&=&
\sum_{j\in I,l\in K_j,\atop{k\in K_j\cap g^{-1}K_j, gk= l}}(*)+
\sum_{j\in I,l\in K_j,\atop{k\in K_j\cap g^{-1}K_j, gk\ne l}}(*)\\
&&+\sum_{j\in I,l\in K_j}\sum_{j\ne i\in I, \atop{k\in K_i\cap g^{-1}K_i}}(*)
     +\sum_{j\in I,l\in K_j}\sum_{i,k\in K_i\backslash g^{-1}K_i}(*) \\
&=&\sum_1(*)+\sum_2(*)+\sum_3(*).
     \end{eqnarray*}
In $\sum_1$ we sum for $i=j, k\in K_i\cap g^{-1}K_i, gk=l$,
in $\sum_2$ we sum for $i=j, k\in K_i\cap g^{-1}K_i, gk\ne l$, or $i\ne j,
k\in K_i\cap g^{-1}K_i, l\in K_j$,
in $\sum_3$ we sum for $k\in K_i\backslash g^{-1}K_i$.
Due to the cocycle identity, we have
$\sum_1(v_g\alpha_g(v_k)v_l^*-1)\alpha_g(E_{i,k})E_{j,l}= 
\sum_{j}(v_g\alpha_g(v_k)v_{gk}^*-1)\alpha_g(E_{i,k})E_{j,l}=0
$, and $\sum_1$ part vanishes. 
Hence 
\begin{eqnarray*}
\lefteqn{|v_g\alpha_g(W^*)W-1|_\phi} \\
&=&|\sum_2(v_g\alpha_g(v_k)v_l^*-1)
\alpha_g(E_{i,k})E_{j,l} 
+\sum_3(v_g\alpha_g(v_k)v_l^*-1)
\alpha_g(E_{i,k})E_{j,l}
|_\phi \\
&\leq& 2\sum_2|\alpha_{g}(E_{i,k})E_{j,l}|_1
+2\sum_3|\alpha_{g}(E_{i,k})E_{j,l}|_1
\end{eqnarray*}
holds by Lemma \ref{lem:ineq}(4).
Similarly we have 
$$ |(v_g\alpha_g(W^*)W-1)^*|_\phi\leq 2\sum_2|\alpha_{g}(E_{i,k})E_{j,l}|_1
+2\sum_3|\alpha_{g}(E_{i,k})E_{j,l}|_1.$$

The estimation of $\sum_3$ is as follows.
\begin{eqnarray*}
 \sum_3|\alpha_g(E_{i,k})E_{j,l}|_1
&=&\sum_{j\in I,l\in K_j}\sum_{i\in I,\atop{k\in K_i\backslash g^{-1}K_j}}
|\alpha_g(E_{i,k})E_{j,l}|_1\\ 
&=&
\sum_{i\in I, \atop{k\in K_j\backslash g^{-1}K_j}}|E_{i,k}|_1 \\
&\leq &\epsilon'+5\epsilon^{'\frac{1}{2}}\leq 6\epsilon^{'\frac{1}{2}},
\end{eqnarray*}
since we have $|K_j\backslash g^{-1}K_j|\leq \epsilon' |K_j|$.

Next we estimate $\sum_2$ part. Then 
\begin{eqnarray*}
 \sum_2|\alpha_g(E_{i,k})E_{j,l}|_1&=&
\sum_{j\in I,\atop{l\in K_j}}\sum_{i\ne j, \atop{k\in K_i\cap g^{-1}K_i}}
|\alpha_g(E_{i,k})E_{j,l}|_1 +
\sum_{j\in I, l\in K_j, \atop{k\in K_j\cap g^{-1}K_j, gk\ne l}}
|\alpha_g(E_{i,k})E_{j,l}|_1 \\ 
&=&\sum_{i\in I, \atop{k\in K_i\cap g^{-1}K_i}}|
\alpha(E_{i,k})(1-\sum_{l}E_{i,l})|_1 +
\sum_{i\in I, \atop{k\in K_i\cap g^{-1}K_i}}|\alpha_g(E_{i,k})
\sum_{l\in K_i, l\ne gk}E_{i,l}|_1 \\
&=& \sum_{i\in I, k\in K_i\cap g^{-1}K_i}|\alpha_g(E_{i,k})(1-E_{i,gk})|_1\\
&=& \sum_{i\in I, k\in K_i\cap g^{-1}K_i}|(\alpha_g(E_{i,k})-
E_{i,gk})(1-E_{i,gk})|_1\\ 
&\leq& 
\sum_{i\in I, k\in K_i\cap g^{-1}K_i}|(\alpha_g(E_{i,k})-
E_{i,gk})|_1\\
&\leq& 10\epsilon^{'\frac{1}{2}} 
\end{eqnarray*}
holds. 
Then we get $|v_g\alpha_g(W^*)W-1|_{\phi}\leq 32\epsilon^{'\frac{1}{2}}$
and $|(v_g\alpha(W^*)W-1)^*|_{\phi}\leq 32\epsilon^{'\frac{1}{2}}$.
By Lemma \ref{lem:ineq} (3),
we have  
$$\|v_g\alpha_g(W^*)W-1\|_\phi^{\#}\leq8\epsilon^{'\frac{1}{4}}<\epsilon.$$
We choose a representing sequence $W=(w_n)$ consisting of unitaries, and 
$E_{i,k}=(e^n_{i,k})$ such that 
$\{e_{i,k}^n\}_{i\in I, k\in K_i}$ is a partition of unity for each $n$.
Set
$a_n:=\sum_{i,k}v_k^*e^n_{i,k}$. 
(Note that $a_n$ is not necessary a unitary.)
Since $W=(w_n)=(a_n)$ in $M^\omega$,   
$\{w_n-a_n\}$ converges to 0 $\sigma$-strongly$*$. 
Choose a sufficiently large $n$ such that 
$$\|v_g\alpha_g(w^*_{n})w_n-1\|_\phi^{\#}<\epsilon,\, g\in F, \phi\in \Phi^+,$$
$$\|[\phi, e_{i,k}^n]\|<\delta,\, \phi\in \Phi^+, i \in I, k\in K_i,$$
$$\|w_n-a_n\|_\phi^{\#}<\epsilon/3, \, \phi\in \Phi^+.$$
Set $w:=w_n$, $a:=a_n$, $e_{i,k}:=e^n_{i,k}$.
(Note that we never use the approximate
commutativity of $v_g$ with elements in $\Psi$ in the estimate of
$\|v_g\alpha_g(w^*)w-1\|_{\phi}^\#$.) 
Next we estimate $\|[a,\varphi]\|$, $\varphi\in \Phi^+$. To do so,  
we estimate $\|[v_g,\varphi]\|$, $g\in K_i$, $\varphi\in \Phi^+$ at first.
We express $g$ as $g=g_1g_2\cdots g_k$, $g_i\in F$, $k\leq L$.
Then it is easily shown $\|[v_g,\varphi]\|\leq
\sum_{i=1}^k\|[v_{g_i},\alpha_{g_1\cdots
g_{i-1}}^{-1}(\varphi)]\|$. (When $i=1$, $\alpha_{g_1\cdots g_{0}}(\varphi)$
means $\varphi$.)
By the definition of $\Psi$, each $\alpha_{g_1\cdots  
g_i}^{-1}(\varphi)$ is in $\Psi$, and by the assumption on $v_g$, $\|[v_{g_i},
\alpha_{g_1\cdots g_{i-1}}^{-1}(\varphi)]\|<\delta$ follows. Hence we have $\|[v_g,
\varphi]\|\leq L\delta $. 
Finally, 
\begin{eqnarray*}
 \|[a,\varphi]\|&\leq& \sum_{i\in I,k\in K_i}\|[v_k^*e_{i,k},\varphi]\| \\
&\leq& \sum_{i\in I,k\in K_i}\|[v_k^*,\varphi]e_{i,k}\|+\|v_k^*[e_{i,k},\varphi]\| \\
&\leq& \sum_{i\in I,k\in K_i}(L+1)\delta \\
&=&\sum_{i\in I}|K_i|(L+1)\delta \\
&< &\epsilon/3
\end{eqnarray*}
holds for $\varphi\in \Phi^+$. 

By Lemma \ref{lem:ineq} (1),
\begin{eqnarray*}
 \|[w, \varphi]\|&\leq& \|[a,\varphi]\|+ \|[w-a,\varphi]\|\\
&<& \epsilon/3 +2\|w-a\|_\varphi^{\#} \\
&< &\epsilon 
\end{eqnarray*}
holds for $\varphi \in \Phi^+$,
and $w$ is a desired unitary. \hfill$\Box$

\noindent
\textbf{Remark.} 
If we replace $v_g\alpha_g(W^*)W$
with $Wv_g\alpha_g(W^*)$ in the above proof, then we get a similar result as in Lemma
\ref{lem:stab} by replacing $\|v_g\alpha_g(w^*)w-1\|_\phi^\#$ with
$\|wv_g\alpha_g(w^*)-1\|_\phi^\#$.   

Now we can present a proof of  Theorem \ref{thm:main} by means of Lemma
\ref{lem:stab}. 

\noindent
\textbf{Proof of Theorem \ref{thm:main}.}
Fix a faithful normal state $\varphi_0$. 
Let $\Phi=\{\varphi_i\}_{i=0}^\infty$ be a countable dense subset of 
$M_*$, and set $\Phi_n:=\{\varphi_i\}_{i=0}^n$. 
Let $G_n \subset\subset G$ be such that $G_n\subset
G_{n+1}$, $G_n^{-1}=G_n$ and $\cup_n G_n=G$.

We construct $w_n, v^n_g, \bar{v}^n_g \in U(M)$, 
$\Phi_n',\Psi_n\subset\subset  M_*$, $\Phi^+_n\subset \subset M_*^+$,   
$\delta_n>0$, actions $\alpha^{(2n)}_g, \beta^{(2n-1)}_g$ of $G$
satisfying the below conditions inductively.
(We set $\alpha^{(0)}:=\alpha$, $\beta^{(-1)}:=\beta$.) 
\begin{eqnarray*}
&&(1.2n)\,  \|\beta^{(2n-1)}_g(\varphi)-\alpha^{(2n)}_g(\varphi)\|<1/2^{2n}, 
g\in G_{2n},\varphi\in \Phi'_{2n}, (n\geq 1). \\ 
&&(1.2n+1)\, \|\alpha^{(2n)}_g(\varphi)-\beta^{(2n+1)}_g
(\varphi)\|<1/2^{2n+1}, 
g\in G_{2n+1},\varphi\in \Phi'_{2n+1}, (n\geq 0). \\
&&(2.2n)\,\|\beta^{(2n-1)}_g(\psi)-\alpha^{(2n)}_g(\psi)\|
<\frac{\delta_{2n-1}}{2}, 
g\in G_{2n-1},
\psi\in \bigcup_{g\in G_{2n-1}}\beta^{(2n-1)}_{g^{-1}}(\Psi_{2n-1}), (n\geq 1). \\
&&(2.2n+1)\,\|\alpha^{(2n)}_g(\psi)-\beta^{(2n+1)}_g(\psi)\|
<\frac{\delta_{2n}}{2},\,g\in G_{2n},
\psi\in \bigcup_{g\in G_{2n}}\alpha^{(2n)}_{g^{-1}}(\Psi_{2n}), (n\geq 1). \\
&&(3.n)\, \|v_g^{n}-1\|_{\phi}^{\#}<1/4^n, g\in G_{n-2},  
\phi\in \Phi_{n-2}^+, (G_{-1}=G_0=G_1, \Phi_{-1}^+=\Phi_0^+=\{\varphi_0\}),
\\
&&(4.n)\, \|[w_n,\varphi]\|<1/4^n, \varphi\in \Phi'_{n-1}, (n\geq 3).\\
&&\bar{v}^{n}_g:=v_{g}^n\Ad w_{n}(\bar{v}^{n-2}_g), 
(\bar{v}^1_g=v_g^1, \bar{v}^2_g=v^2_g), \\
&&\alpha_g^{(2n)}:=\Ad v^{2n}_g\circ\Ad w_{2n}^*\circ\alpha_g^{(2n-2)}\circ 
 \Ad w_{2n},\\
&& \beta_g^{(2n-1)}:=\Ad {v}^{2n-1}_g\circ\Ad w_{2n-1}^*\circ\beta_g^{(2n-3)}\circ 
 \Ad w_{2n-1},\\
&&\Phi_{2n}':=\Phi_{2n}\cup\Ad w_{2n-1}^*w_{2n-3}^*\cdots
w_1^*(\Phi_{2n})\cup\{\bar{v}^{2n-1}_g\varphi_0,
\varphi_0\bar{v}^{2n-1}_g\}_{g\in G_{2n-1}}, \\
&&\Phi_{2n+1}':=\Phi_{2n+1}\cup\Ad w_{2n}^*w_{2n-2}^*\cdots
w_2^*(\Phi_{2n+1})\cup\{\bar{v}^{2n}_g
\varphi_0,\varphi_0\bar{v}^{2n}_g\}_{g\in G_{2n}}, (n\geq 1) \\
&&\Phi_n^+:=\{\Ad \bar{v}_g^{n}(\varphi_0)\mid g\in G_n\}.
\end{eqnarray*}

Here $\delta_{2n}$ and $\Psi_{2n}$ ($\delta_{2n-1}$ and $\Psi_{2n-1}$)
are chosen as in Lemma
\ref{lem:stab} for $\alpha^{(2n)}$, 
$1/4^{2n+2}$,
$G_{2n}$, $\Phi_{2n}^+$ and
$\Phi'_{2n+1}$ (resp. for $\beta^{(2n-1)}$, $1/4^{2n+1}$, $G_{2n-1}$,
$\Phi_{2n-1}^+$, and 
$\Phi'_{2n}$).

First set $\Phi'_1:=\Phi_1$ and fix a $\beta$-cocycle $u^1_g$ such that 
$$\|\alpha^{(0)}_g(\varphi)-\Ad u_g^1\beta^{(-1)}_g(\varphi)\|< 1/2,\,
g\in G_1,  
\varphi\in \Phi'_1.$$
By Lemma \ref{lem:stab}, we get a unitary $w_1$ such that
$\|u_g\beta^{(-1)}_g(w_1^*)w_1-1\|_{\varphi_0}^{\#}<1/4$, $g\in G_1$. Set
$v^1_g:=u^1_g\beta_g^{(-1)}(w_1^*)w_1 $, and
$\beta_g^{(1)}:=\Ad u^1_g\circ\beta_g^{(-1)}=\Ad 
v_g^1 \circ
\Ad w_1^*\circ\beta_g^{(-1)}\circ \Ad w_1$. 
Then we have 
$$(1.1)\quad
\|\alpha_g^{(0)}(\varphi)-\beta_g^{(1)}(\varphi)\|<1/2$$ and
$$(3.1)\quad \|v^1_g-1\|_{\varphi_0}^{\#}<1 /4$$
 for $g\in G_1$.  
Set $\bar{v}_g^1:=v^1_g$, 
$\Phi_2':=\Phi_2\cup \Ad w_1^*(F_2) \cup \{\bar{v}_g^1\varphi_0,
\varphi_0\bar{v}_g^1\}$, and $\Phi_1^+:=
\{\Ad \bar{v}_g^1(\varphi_0)\}_{g\in G_1}$.
By Lemma \ref{lem:stab}, we choose $\Psi_1$ and $\delta_1$ for
$\beta_g^{(1)}$, $1/4^3$, $G_1$, $\Phi_1^+$ and $\Phi'_2$.

Next we take an $\alpha^{(0)}$-cocycle $u^2_g$ such that 
\begin{eqnarray*}
(a.2)&& \|\beta_g^{(1)}(\varphi)-\Ad u^2_g\alpha_g^{(0)}(\varphi)\|
<\frac{1}{2^2},  
\,g\in G_2, \varphi\in \Phi'_2, \\
(b.2)&& \|\beta_g^{(1)}(\psi)-\Ad u^2_g\alpha_g^{(0)}(\psi)\|
<\frac{\delta_1}{2}, \,g\in G_1,\psi\in
\bigcup_{g\in G_1}\beta_{g^{-1}}^{(1)}(\Psi_1).  
\end{eqnarray*}

By Lemma \ref{lem:stab}, we get $w_2\in U(M)$ such that 
$\|u_g^2\alpha_0(w_2^*)w_2-1\|_{\varphi_0}^{\#}<1/4^2$ for
$g\in G_2$. Set
$v_g^2=\bar{v}_g^2:=u_g^2\alpha_g^{(0)}(w_2^* )w_2$ 
and $\alpha_g^{(2)}:=\Ad u_g^2\alpha_g^{(0)}
=\Ad v_g^2\circ \Ad w_2^*\circ \alpha_g^{(0)}\circ \Ad w_2$.  
Then  we get
$$(3.2) \quad\|v_g^2-1\|_{\varphi_0}^{\#}<\frac{1}{4^2},\, 
g\in G_2.$$\\  
By $(a.2)$ and $(b.2)$,  
\begin{eqnarray*}
(1.2)&&\quad\|\beta_g^{(1)}(\psi)-\alpha_g^{(2)}(\psi)\|<\frac{1}{2^2}, \,
g\in G_2,\varphi\in \Phi'_2, \\ 
(2.2)&&\quad\|\beta_g^{(1)}(\psi)-\alpha_g^{(2)}(\psi)\|<
\frac{\delta_1}{2},\,g\in G_1, \psi\in\bigcup_{g\in G_1}
\beta_{g^{-1}}^{(1)}(\Psi_1). 
\end{eqnarray*}

Set $\Phi'_3:=\Phi_3\cup \Ad w_2^*(\Phi_3) \cup
\{\bar{v}_g^2\varphi_0,\varphi_0\bar{v}_g^2\}_{g\in G_2}$ and 
$\Phi_2^+:=\{\Ad \bar{v}_g^2 \varphi_0\mid g\in G_2 \}$. 
By Lemma \ref{lem:stab}, 
we choose $\delta_2$ and $\Psi_{2}$ for $\alpha^{(2)}$, $G_2$,
$\Phi'_3$, $\Phi_2^+$ and $1/4^4$.

Suppose that we construct up to $\alpha_g^{(2n)}$, $\beta_g^{(2n-1)}$,
$w_{2n}$, $v_g^{2n}$, $\bar{v}_g^{2n}$, 
$\Phi'_{2n+1}$, $\delta_{2n}$ and $\Psi_{2n}$.

We choose a $\beta^{(2n-1)}$-cocycle $u^{2n+1}_g\in U(M)$ such that 
\begin{eqnarray*}
(a.2n+1)&& \|\alpha_g^{(2n)}(\varphi)-\Ad
u_g^{2n+1}\beta_g^{(2n-1)}(\varphi)\|<\frac{1}{2^{2n+1}},\,g\in G_{2n+1},
\varphi\in \Phi'_{2n+1},\\ 
(b.2n+1)&& \|\alpha_g^{(2n)}(\psi)-\Ad
u_g^{2n+1}\beta_g^{(2n-1)}(\psi)\|<\frac{\delta_{2n}}{2},\,g\in G_{2n},
 \psi\in 
\bigcup_{g\in G_{2n}}\alpha_{g^{-1}}^{(2n)}(\Psi_{2n}),\\ 
(c.2n+1)&& \|\alpha_g^{(2n)}(\psi)-\Ad
u_g^{2n+1}\beta_g^{(2n-1)}(\psi)\|<\frac{\delta_{2n-1}}{2},\,g\in G_{2n-1},
\varphi\in \bigcup_{g\in G_{2n-1}}\beta_{g^{-1}}^{(2n-1)}(\Psi_{2n-1}).
\end{eqnarray*}

Then by $(2.2n)$ and $(c.2n+1)$ we get 
$$\|\beta_g^{(2n-1)}(\psi)-\Ad
u_g^{2n+1}\beta_g^{(2n-1) }(\psi)\|<\delta_{2n-1}, \,g\in
G_{2n-1},\psi\in 
\beta^{(2n-1)}_{g^{-1}}(\Psi_{2n-1}),$$ which yields
$\|[u_g^{2n+1},\psi]\|<\delta_{2n-1}, g\in G_{2n-1},\psi\in
\Psi_{2n-1}$. By the
choice of $\delta_{2n-1}$ and $\Psi_{2n-1}$, there exists a unitary
$w_{2n+1}$ such that
$$\|u_g^{2n+1}\beta_g^{(2n-1)}(w_{2n+1}^*)w_{2n+1}-1\|_{\phi}^{\#}
<1/4^{2n+1},\,g\in G_{2n-1},\phi\in \Phi_{2n-1}^+,$$  
and $$(4.2n+1)\quad\|[w_{2n+1}, 
\varphi]\|<1 /4^{2n+1},\varphi\in \Phi'_{2n}.$$  Set
$v_g^{2n+1}:=u_g^{2n+1}\beta_g^{(2n-1)}( w_{2n+1}^*)w_{2n+1}$ and
$\beta_g^{(2n+1)}:=\Ad u_g^{(2n+1)}\beta_g^{(2n-1)}
=\Ad v_g^{(2n+1)}\circ\Ad w_{2n+1}^*\circ\beta_g^{(2n-1)}\circ\Ad
w_{2n+1}$.  Then
$$(3.2n+1)\quad\|v_g^{2n+1}-1\|_{\phi}<\frac{1}{4^{2n+1}},\, 
g\in G_{2n-1},\phi\in \Phi_{2n-1}^+.$$
By $(a.2n+1)$
and $(b.2n+1)$, we get 
\begin{eqnarray*}
(1.2n+1)&&\|\alpha_g^{(2n)}(\varphi)-\beta_g^{(2n+1)}(\varphi)\|
<\frac{1}{2^{2n+1}},\,g\in G_{2n+1}, 
\varphi\in \Phi'_{2n+1}. \\ 
(2.2n+1)&&\|\alpha_g^{(2n)}(\psi)-\beta_g^{(2n+1)}(\psi)\|
<\frac{\delta_{2n}}{2}, \, g\in G_{2n},
\psi\in\bigcup_{g\in G_{2n}}\alpha_{g^{-1}}^{(2n)}(\Psi_{2n}).  
\end{eqnarray*} 

Set $\bar{v}_g^{2n+1}:=v_g^{2n+1}\Ad w_{2n+1}^*(\bar{v}_g^{2n-1})$, and
define  
$$\Phi'_{2n+2}:=\Phi_{2n+2}\cup \Ad w_{2n+1}^*w_{2n-1}^*\cdots
w_1^*(\Phi_{2n+2})\cup\{\bar{v}_g^{2n+1}(\varphi_0),
\varphi_0\bar{v}_g^{2n+1}\}_{g\in G_{2n+1}},$$
$$\Phi_{2n+1}^+:=\{\Ad \bar{v}_g^{2n+1}\varphi_0\mid g\in G_{2n+1}\}.$$   

By Lemma \ref{lem:stab}, we choose $\delta_{2n+1}>0$ and
$\Psi_{2n+1}\subset\subset M_*$ for $\beta^{(2n+1)}$,
$1/4^{2n+3}$, $G_{2n+1}$, $\Phi_{2n+1}^+$
and $\Phi' _{2n+2}$, and the $2n$-th step is finished. 

Next we choose an $\alpha^{(2n)}$-cocycle $u_g^{2n+2}$ such that 
\begin{eqnarray*}
&&(a.2n+2)\, \|\beta_g^{(2n+1)}(\varphi)-\Ad u_g^{2n+2}
\alpha_g^{(2n)}(\varphi)\|<\frac{1}{2^{2n}},\,
g\in G_{2n+2},\varphi\in \Phi'_{2n+2}, \\
&&(b.2n+2)\, \|\beta_g^{(2n+1)}(\psi)-\Ad u_g^{2n+2}
\alpha_g^{(2n)}(\psi)\|<\frac{\delta_{2n+1}}{2},\,g\in G_{2n+1},
\psi\in\bigcup_{g\in G_{2n+1}} \beta^{(2n+1)}_{g^{-1}}(\Psi_{2n+1}), \\
&&(c.2n+2)\, \|\beta_g^{(2n+1)}(\psi)-
\Ad u_g^{2n+2}\alpha_g^{(2n)}(\psi)\|<\frac{\delta_{2n}}{2},\,g\in G_{2n},
\psi\in \bigcup_{g\in G_{2n}}\alpha_{g^{-1}}^{(2n)}(\Psi_{2n}). 
\end{eqnarray*}

By $(c.2n+2)$ and $(2.2n+1)$, we get 
$$\|\alpha_g^{(2n)}(\psi)-\Ad
u_g^{2n+2}\alpha_g^{(2n)}(\psi)\|<\delta_{2n},\,g\in
G_{2n},\psi\in \bigcup_{g\in G_{2n}}\alpha_{g^{-1}}^{(2n)}(\Psi_{2n}).$$
We thus have $\|[u_g^{2n+2},\psi]\|<\delta_{2n}$ for $g\in G_{2n}$
and $\psi\in \Psi_{2n}$. By the choice of $\delta_{2n}$ and $\Psi_{2n}$, we can find 
$w_{2n+2}\in U(M)$ such that
$$\|u_g^{2n+2}\alpha_g^{(2n)}(w_{2n+2}^*)w_{2n+2}-1\|^\#_{\phi}<
\frac{1}{4^{2n+2}},\,g\in G_{2n},\phi\in \Phi_{2n}^+$$
and $$(4.2n+2)\quad \|[w_{2n+2},\varphi]\|<\frac{1}{4^{2n+2}},\, 
\varphi\in \Phi'_{2n+1}.$$
Set ${v}_g^{2n+2}:=u_g^{2n+2}\alpha_g^{(2n)}(w_{2n+2}^*)w_{2n+2}$ and
$\alpha_g^{(2n+2)}:=\Ad u_g^{2n+2}\alpha_g^{(2n)}=\Ad v_g^{2n+2}\circ
\Ad w_{2n+2}^*\circ \alpha_g^{(2n)}\circ\Ad w_{2n+2} $.
Then $$(3.2n+2)\quad \|v_g^{2n+2}-1\|^\#_{\phi}<
\frac{1}{4^{2n+2}},\,g\in
G_{2n}, \phi\in \Phi^+_{2n},$$ and 
by $(a.2n+2)$ and $(b.2n+2)$, we get
\begin{eqnarray*}
(1.2n+2)&&\|\beta_g^{(2n+1)}(\varphi)-\alpha_g^{(2n+2)}(\varphi)\|<
\frac{1}{2^{2n}},
\,g\in G_{2n},\varphi\in \Phi'_{2n+2}. \\
(2.2n+2)&& \|\beta_g^{(2n+1)}(\psi)-\alpha_g^{(2n+2)}(\psi)\|
<\frac{\delta_{2n+1}}{2},
\,g\in G_{2n+1},\psi\in\bigcup_{g\in G_{2n+1}} 
\beta_{g^{-1}}^{(2n+1)}(\Psi_{2n+1}). 
\end{eqnarray*}

Set $\bar{v}_g^{2n+2}:=v_g^{2n+2}\Ad w_{2n+2}^*(\bar{v}_g^{2n})$ and 
$$\Phi'_{2n+3}:=\Phi_{2n+3}\cup \Ad w_{2n+2}^*w_{2n}^*\cdots
w_2^*(\Phi_{2n+3})\cup\{\bar{v}_g^{2n+2}\varphi_0,
\varphi_0\bar{v}_g^{2n+2}\},$$
$$\Phi_{2n+2}^+=\{\Ad \bar{v}_g^{2n+2}\varphi_0\mid g\in G_{2n+2} \}.$$ 

We choose $\delta_{2n+2}>0$, $\Psi_{2n+2}\subset\subset M_*$ for $\alpha_{g}^{(2n)}$,
$1/4^{2n+4}$, $G_{2n+2}$, 
$\Phi^+_{2n+2}$ and  $\Phi'_{2n+3}$ by Lemma \ref{lem:stab}. Then the
$(2n+1)$-th  step is finished, 
and thus we complete induction.

Set $\theta_{2n}:=\Ad w_{2n}^*w_{2n-2}^*\cdots w_2^*$. Then 
we have $\alpha_g^{(2n)}=\Ad \bar{v}_g^{2n}\circ\theta_{2n}\circ
\alpha_g\circ \theta_{2n}^{-1}$.  
We will verify $\{\theta_{2n}\}$ converges to some $\theta\in \Aut(M)$.
To do so, we will prove $\{\theta_{2n}(\varphi)\}$ and $\{\theta_{2n}^{-1}(\varphi)\}$
are Cauchy sequences for $\varphi\in M_*$.
Suppose $\varphi\in \Phi_k$. For any $n$ with $k\leq 2n+1$, $\varphi$ and
$\theta_{2n}(\varphi)$ are in $ 
\Phi'_{2n+1}$. 
By $(4.2n+2)$, we have
\begin{eqnarray*}
 \|\theta_{2n+2}(\varphi)-\theta_{2n}(\varphi)\|&=&\|[w_{2n+2},
\theta_{2n}(\varphi)]\| \\
&<&\frac{1}{4^{2n+2}}, 
\end{eqnarray*}
and
\begin{eqnarray*}
 \|\theta_{2n+2}^{-1}(\varphi)-\theta_{2n}^{-1}(\varphi)\|&=&
\|w_{2n+2}^*\varphi w_{2n+2}-\varphi\| \\
&< &\frac{1}{4^{2n+2}}. 
\end{eqnarray*}
It follows that $\{\theta_{2n}(\varphi)\}$ and $\{\theta_{2n}^{-1}(\varphi)\}$ 
are Cauchy sequences for $\varphi\in \Phi$. Then so are 
$\{\theta_{2n}(\varphi)\}$ and $\{\theta_{2n}^{-1}(\varphi)\}$
for every $\varphi\in M_*$, 
since $\Phi$ is dense in $M_*$.
Hence $\{\theta_{2n}\}$ converges to 
some $\theta\in \Aut(M)$.

Next we will verify $\{\bar{v}_g^{2n}\}$ is a Cauchy sequence with respect to
$\|\cdot\|_{\varphi_0}^{\#}$.
Since we have
\begin{eqnarray*}
\|\bar{v}_g^{2n+2}-\bar{v}_g^{2n}\|_{\varphi_0}^{\#} 
&\leq&\|(v_g^{2n+2}-1)\bar{v}_g^{2n}\|_{\varphi_0}^{\#} 
+ \|(w_{2n+2}\bar{v}_g^{2n}w_{2n+2}^*-\bar{v}_g^{2n})\|_{\varphi_0}^{\#}
 \\
&&+\|(v_g^{2n+2}-1)(w_{2n+2}\bar{v}_g^{2n}w_{2n+2}^*-\bar{v}_g^{2n})
\|_{\varphi_0}^{\#}, 
\end{eqnarray*}
we will estimate the above three terms. 

Suppose $g\in G_{k}$. Then for any $n$ with $2n\geq k$,
$\varphi_0,\Ad \bar{v}_{g}^{2n}(\varphi_0)\in
\Phi^+_{2n}$, and hence $\|v_g^{2n+2}-1\|^\#_{\Ad \bar{v}_g^{2n}(\varphi_0)}<
1 /4^{2n+2}$ and
$\|v_g^{2n+2}-1\|^\#_{\varphi_0}<
1 /4^{2n+2}$ hold
by $(3.2n+2)$.  

We have
\begin{eqnarray*}
 \|(v_g^{2n+2}-1)\bar{v}_g^{2n}\|_{\varphi_0}^{\#2} &=& 
\frac{1}{2}(\|(v_g^{2n+2}-1)\bar{v}_g^{2n}\|_{\varphi_0}^{2}+
\|\bar{v}_g^{2n*}(v_g^{2n+2*}-1)\|_{\varphi_0}^{2} ) \\ 
&= &\frac{1}{2}(\|(v_g^{2n+2}-1)\|_{\Ad \bar{v}_g^{2n}(\varphi_0)}^{2}+
\|(v_g^{2n+2*}-1)\|_{\varphi_0}^{2} ) \\ 
&\leq & \|(v_g^{2n+2}-1)\|_{\Ad \bar{v}_g^{2n}(\varphi_0)}^{\#2}+
\|(v_g^{2n+2}-1)\|_{\varphi_0}^{\#2}  \\
&< &\frac{2}{16^{2n+2}}. 
\end{eqnarray*}
Hence we get $\|(v_g^{2n+2}-1)\bar{v}_g^{2n}\|_{\varphi_0}^{\#}
<\sqrt{2}/4^{2n+2}<1/2^{2n+1}.  $

We next estimate 
$\|w_{2n+2}^*\bar{v}_g^{2n}w_{2n+2}-\bar{v}_{2n}\|_{\varphi_0}^{\# }$.
Since $\varphi_0,\bar{v}_g^{2n}\varphi_0, \varphi_0\bar{v}_g^{2n}\in
\Phi'_{2n+1}$, we have $\|[w_{2n+2},\varphi_0]\|<1/4^{2n+2}$, 
$\|[w_{2n+2},\bar{v}_g^{2n}\varphi_0]\|<1/4^{2n+2}$ and 
$\|[w_{2n+2},\varphi_0\bar{v}_g^{2n}]\|<1/4^{2n+2}$ by
$(4.2n+2)$.
Then 
\begin{eqnarray*}
\lefteqn{
\|(w_{2n+2}^*\bar{v}_g^{2n}w_{2n+2}-\bar{v}_g^{2n})\cdot\varphi_0 \|} \\&=&
 \|(w_{2n+2}^*\bar{v}_g^{2n}-\bar{v}_g^{2n}w_{2n+2}^*)w_{2n+2}\cdot\varphi_0 \|\\
 &\leq&
\|(w_{2n+2}^*\bar{v}_g^{2n}-\bar{v}_g^{2n}w_{2n+2}^*)\cdot\varphi_0\cdot w_{2n+2} \|+
\|(w_{2n+2}^*\bar{v}_g^{2n}-\bar{v}_g^{2n}w_{2n+2}^*)\cdot[w_{2n+2},\varphi_0] \|
\\ &\leq& 
\|w_{2n+2}^*\bar{v}_g^{2n}\varphi_0w_{2n+2}-\bar{v}_g^{2n}w_{2n+2}^*
\varphi_0w_{2n+2} \|+
 2/4^{2n+2} \\
&\leq &
\|[w_{2n+2}^*,\bar{v}_g^{2n}\varphi_0]w_{2n+2}\|+
\|\bar{v}_g^{2n}\varphi_0-\bar{v}_g^{2n}w_{2n+2}^*
\varphi_0w_{2n+2} \|+
 2/4^{2n+2} \\
&\leq &
1/4^{2n+2}+
\|\bar{v}_g^{2n}[\varphi_0,w_{2n+2}^*] \|+
 2/4^{2n+2} \\
&\leq &1/4^{2n+1}
\end{eqnarray*}
holds. Hence 
\begin{eqnarray*}
 \|w_{2n+2}^*\bar{v}_g^{2n}w_{2n+2}-\bar{v}_g^{2n}\|^2_{\varphi_0}&\leq&
\| w_{2n+2}^*\bar{v}_g^{2n}w_{2n+2}-\bar{v}_g^{2n}\| 
\|(w_{2n+2}^*\bar{v}_g^{2n}w_{2n+2}-\bar{v}_g^{2n})\cdot\varphi_0\|\\
&\leq&
2\|(w_{2n+2}^*\bar{v}_g^{2n}w_{2n+2}-\bar{v}_g^{2n})\cdot\varphi_0\|\\
&\leq & \frac{2}{4^{2n+1}}
\end{eqnarray*}
holds.

In a similar way, we can show 
$\|(w_{2n+2}^*\bar{v}_g^{2n}w_{2n+2}-\bar{v}_g^{2n})^*\|^2_{\varphi_0}\leq 
2/4^{2n+1}$.
Hence we get
$\|w_{2n+2}^*\bar{v}_g^{2n}w_{2n+2}-\bar{v}_g^{2n}\|^\#_{\varphi_0}\leq
\sqrt{2/4^{2n+1}}=\sqrt{2}/{2^{2n+1}}<1/
2^{2n}$.

The third term  
$\|(v_g^{2n+2}-1)(w_{2n+2}^*\bar{v}_g^{2n}w_{2n+2}-
\bar{v}_g^{2n})\|_{\varphi_0}^{\#}$ is estimated as follows.  
\begin{eqnarray*}
 \|(v_g^{2n+2}-1)(w_{2n+2}^*\bar{v}_g^{2n}w_{2n+2}-
\bar{v}_g^{2n})\|_{\varphi_0}^{\# 2}
&=&\frac{1}{2}(
\|(v_g^{2n+2}-1)(w_{2n+2}^*\bar{v}_g^{2n}w_{2n+2}-
\bar{v}_g^{2n})\|_{\varphi_0}^{2} \\&&+
\|(w_{2n+2}^*\bar{v}_g^{2n*}w_{2n+2}-
\bar{v}_g^{2n*})(v_g^{2n+2*}-1)\|_{\varphi_0}^{2} ) \\
&\leq & 2\|w_{2n+2}^*\bar{v}_g^{2n}w_{2n+2}-\bar{v}_g^{2n}\|_{\varphi_0}^{2}+
2\|v_g^{2n+2*}-1\|_{\varphi_0}^2 \\
&\leq& 4\|w_{2n+2}^*\bar{v}_g^{2n}w_{2n+2}-\bar{v}_g^{2n}\|_{\varphi_0}^{\#2}+
4\|v_g^{2n+2*}-1\|_{\varphi_0}^{\#2} \\
&\leq& \frac{4}{4^{2n}}+\frac{4}{4^{2n+2}} \\
&\leq & \frac{2}{4^{2n-1}}.
\end{eqnarray*}
Hence $\|(v_g^{2n+2}-1)(w_{2n+2}^*\bar{v}_g^{2n}w_{2n+2}-
\bar{v}_g^{2n})\|_{\varphi_0}^{\#}\leq
\sqrt{2}/2^{2n-1}<1/2^{2n-2}$, and 
we have the following.

\begin{eqnarray*}
 \|\bar{v}_g^{2n+2}-\bar{v}_g^{2n}\|_{\varphi_0}^{\#}&\leq & 
\frac{1}{2^{2n+1}}+
\frac{1}{2^{2n}}+\frac{1}{2^{2n-2}} \\
&\leq &\frac{1}{2^{2n-3}}.
\end{eqnarray*}
It follows that  $\{\bar{v}_g^{2n}\}$ is a Cauchy sequence and converges to some unitary
$\hat{v}_g^{0}$.  

In the same way, 
we can show $\sigma_{2n+1}:=\Ad w_{2n+1}^*w_{2n-1}^*\cdots w_1^*$ and
$\bar{v}_g^{2n+1}$ converges to $\sigma \in \Aut (M)$ and $\hat{v}_g^{1}\in
U(M)$ respectively.  By $(1.n)$ we get 
$\Ad \hat{v}_g^{0}\circ\theta\circ\alpha_g\circ\theta^{-1}=
\Ad \hat{v}_g^{1}\circ \sigma\circ\beta_g\circ \sigma^{-1}$, and hence  
$\alpha$ and $\beta$ are cocycle conjugate. By the construction,
$\theta,\sigma$ are approximately inner.

We will choose a cocycle close to 1. Suppose $\Ad v_g\alpha=\theta\circ
\beta_g\circ \theta^{-1}$, $\theta\in \overline{\Int}(M)$. Fix $F\subset\subset G$. 
Then there exists a unitary $w$ such that
$\|wv_g\alpha_g(w^*)-1\|_{\varphi_0}^\#<\varepsilon$ for each $g\in F$. 
(See the remark after Lemma \ref{lem:stab}.) 
Define a new $\alpha$-cocycle $v_g'$ by
$v'_g:=wv_g\alpha_g(w^*)$. 
We then have $\|v_g'-1\|_{\varphi_0}^\#<\varepsilon$ for $g\in F$, and
\begin{eqnarray*}
\Ad v_g'\alpha_g&=&\Ad(wv_g\alpha_g(w^*))\circ\alpha_g \\
&=&\Ad w\circ \Ad v_g\circ\alpha_g\circ \Ad w^*  \\
&=& \Ad w\circ\theta\circ\beta_g\circ \theta^{-1}\circ\Ad w^*.
\end{eqnarray*} 
Put $\sigma:=\Ad w\circ\theta$. Then $\sigma\in \overline{\Int}(M)$, and 
we get $\Ad
v_g'\alpha_g=\sigma\circ\beta_g\circ\sigma^{-1}$. 
\hfill$\Box$

We present applications of Theorem \ref{thm:main}.
Let $M$ be an injective factor. By the Connes-Krieger-Haagerup
classification of injective factors \cite{Co-inj}, \cite{Kri-erg},
\cite{Co-III1}, \cite{Ha-III1}, $M$ 
is a McDuff factor. 
Since $\Ker\,\mod=\overline{\mathrm{Int}}(M)$
by \cite{Con-surv} and \cite{KwST}, we get the following corollary.

\begin{cor}
 Let $M$ be an injective factor, $G$ a discrete amenable group,
 $\alpha,\beta$  centrally free
 actions of $G$ on $M$.
 Then $\Ad v_g\alpha_g=\theta\circ\beta_g\circ\theta^{-1}$ for some
 $\alpha$-cocycle  $v_g$ and $\theta\in \overline{\Int}(M)$
 if and only if $\mod(\alpha)=\mod(\beta)$.  $($In the type II$_1$ case,
 we regard $\mod(\alpha)$ is trivial for $\alpha\in \Aut(M)$.$)$
\end{cor}

Theorem \ref{thm:main} can be modified for a relatively McDuff subfactor
$N\subset M$ by appropriate changes. 
Indeed  in the proof
we  only have to replace $M_\omega$ with 
$M_\omega\cap N^\omega$, 
which is the subfactor version of a central sequence algebra.
Especially if $N\subset M$ is a strongly
amenable subfactor of type II$_1$ in the sense of Popa \cite{Po-amen},
then it is relatively McDuff thanks to Popa's classification theorem of
strongly amenable subfactors of type II$_1$ \cite{Po-amen}.
We also have $\overline{\mathrm{Int}}(M,N)=\mathrm{Ker}\,\Phi$ by
\cite{Loi-auto}, where $\Phi(\alpha)$ is the Loi 
invariant for $\alpha\in \Aut(M,N)$, and the equivalence between strong
outerness and central freeness by \cite{Po-act}. 
(Also see \cite{M-III1} for the latter fact.)
Hence Theorem \ref{thm:main} gives an alternative proof
of the main theorem in \cite{Po-act}.

\begin{cor}\label{cor:sub1}
 Let $N\subset M$ be a strongly amenable subfactor of type II$_1$, $G$
 a discrete amenable group, and
 $\alpha, \beta $ strongly outer actions of $G$ on $N\subset M$. 
Then $\Ad v_g\alpha_g=\theta\circ\beta_g\circ\theta^{-1}$ for some
 $\alpha$-cocycle  $v_g\in U(N)$ and $\theta\in \overline{\Int}(M, N)$ if and  
 only  if $\Phi(\alpha)=\Phi(\beta)$.
\end{cor}  

When $N\subset M$ is a strongly amenable subfactor of type II$_\infty$,
then we have $\overline{\mathrm{Int}}(M,N)=\mathrm{Ker}\,\Phi\cap
\mathrm{Ker}(\mod)$. 
Hence we have the following corollary.

\begin{cor}\label{cor:sub2}
 Let $N\subset M$ be a strongly amenable subfactor of type II$_\infty$, $G$
 a discrete amenable group, and
 $\alpha, \beta $ strongly outer actions of $G$ on $N\subset M$. 
Then $\Ad v_g\alpha_g=\theta\circ\beta_g\circ\theta^{-1}$ for some
 $\alpha$-cocycle  $v_g\in U(N)$ and $\theta\in \overline{\Int}(M, N)$ if and  
 only  if $\Phi(\alpha)=\Phi(\beta)$ and $\mod(\alpha)=\mod(\beta)$.
\end{cor}  

It is worth noting that Corollary \ref{cor:sub2} implies the
classification of strongly amenable subfactor of type III$_\lambda$,
$0<\lambda<1$,  in \cite{Loi-auto}, \cite{Po-act}. 


\ifx\undefined\bysame
\newcommand{\bysame}{\leavevmode\hbox to3em{\hrulefill}\,}
\fi

\end{document}